\documentclass[10pt,twocolumn,twoside]{IEEEtran}
\usepackage{cite}
\usepackage{amsmath,amssymb,amsfonts}
\usepackage{multicol}
\usepackage{subcaption}
\usepackage{array}
\usepackage{graphicx}
\usepackage{xcolor}
\usepackage{textcomp}
\usepackage{mathtools}
\usepackage{algorithm}
\usepackage{algpseudocode}
\usepackage{subcaption}
\usepackage{booktabs}
\usepackage{amsthm}

\newtheorem{theorem}{Theorem}
\newtheorem{remark}{Remark}
\newtheorem{assumption}{Assumption}
\newtheorem{definition}{Definition}
\newtheorem{lemma}{Lemma}
\newtheorem{corollary}{Corollary}

\usepackage{ragged2e}

\newenvironment{problem}{{\bf Problem. }}{}
\newcommand\Mycomb[2][^n]{\prescript{#1\mkern-0.5mu}{}C_{#2}}
\def\BibTeX{{\rm B\kern-.05em{\sc i\kern-.025em b}\kern-.08em
    T\kern-.1667em\lower.7ex\hbox{E}\kern-.125emX}}
\begin{document}
\title{Minimum time consensus for damped second order agents using Gr\"{o}bner basis}
\author{Akansha Rautela, Deepak U. Patil, Ameer Mulla and Indra Narayan Kar
\thanks{A.Rautela, D.U.Patil and I.N.Kar are with the Department of Electrical Engineering, Indian Institute of Technology Delhi, New Delhi, India.}      
\thanks{A.Mulla is with Department of Electrical Engineering, Indian Institute of Technology Dharwad, Karnataka, India.}
\thanks{email: Akansha.Rautela@ee.iitd.ac.in, deepakpatil@ee.iitd.ac.in, ameer@iitdh.ac.in, ink@ee.iitd.ac.in}}

\maketitle

\begin{abstract}
A problem of achieving minimum time consensus for a set of $N$ second-order LTI system agents with bounded inputs and fuel constraints is considered. Unlike most other works, here the damping effect in agent dynamics is included. First, the attainable set for each agent with fuel budget constraints is characterized, and its boundary equations are derived. Then, using the convexity property, the minimum time at which attainable sets of all agents have a non-empty intersection is computed. By applying Helly's theorem, the computation reduces to finding the minimum time to consensus and the corresponding consensus point for each of the $\Mycomb[N]{3}$ triplets separately.
\end{abstract}


\section{Introduction}
\label{sec:introduction}

   Control of a multi-agent system for various goals such as consensus, synchronization, formation, etc. has been studied extensively in the literature \cite{4118472,ren2010distributed,inputconstraintsieeetac2018,safeconsensusbarrierfunctionieeetac2023,energystateconstraintsautomatica2024,flockingieeetac2023}.    
In consensus problems, a control law is sought such that each agent aims to match its state with that of the neighboring agent to reach a common agreement or point \cite{lin2003multi,dong2013leader}. This article considers the problem of computing a consensus point that a \emph{class of second-order agents with damping} can reach in minimum time, subject to \emph{bounded control inputs} and a \emph{limited fuel budget}. Limitations on fuel budget bring in several restrictions with respect to agent initial states and thereby pose several difficulties in achieving consensus. We characterize a region of initial conditions from where consensus can occur and also a region of consensus points. 
 
 The speed of achieving consensus improves with increasing graph connectivity (the second smallest real eigenvalue of the Laplacian matrix, known as the Fielder eigenvalue) \cite{ren2010distributed}. However, constraints on system resources, such as operational time restrictions and fuel or energy usage, lead to limitations on input \cite{minimumswitchthruster,attitudecontrol,fuelconstraintsinUAVs}, which must be taken into account. Recently, input, state, and energy constraints have been considered in consensus control laws in \cite{inputconstraintsieeetac2018}, \cite{energystateconstraintsautomatica2024}, using model transformations and LMI-based approaches. However, optimal control objectives such as \emph{time} and \emph{fuel} have been recently explored in our previous articles \cite{doi:10.1080/00207179.2017.1285054, patil2019distributed,tacsubmitted2024,minfuelardpakink,sclsubmitted2025} for double integrator agents. Our approach in these articles is set-based, wherein we characterize attainable set boundaries and compute the intersection of the attainable set of all agents. In this letter, a similar set-based approach is used to compute the minimum time to consensus and the corresponding consensus time for second-order agents with a damping term.   
 
 Usually, in literature and also in our previous articles, the double integrator model is used due to its simplicity. However, it lacks a damping term that typically appears in situations where drag or friction is present. In this letter, we extend the results presented in \cite{tacsubmitted2024} to the case where a damping term is present in the agent model. The presence of damping terms leads to difficulties in computing the intersection of attainable sets. We need a solution to boundary equations of the attainable set, which are described by equations consisting of exponential terms when a damping term is present. We make use of a simple substitution for converting exponential terms to polynomials and use Gr\"{o}bner basis to solve these boundary equations systematically. 
 

Three main contributions of this work are --  1) a procedure for computing minimum time to consensus and corresponding consensus state for a set of $N$ identical class of damped second order agents each having a fixed fuel budget is developed using Gr\"{o}bner basis and Helly's theorem, 2) the attainable set is algebraically characterized in terms of its boundary equations and inequalities, and 3) Initial conditions for agents from where consensus is possible and the region in the state space where consensus occurs is also characterized.

\section{Problem formulation and preliminaries}
\subsection{Problem Formulation}
Consider a set of $N$ agents described by the second-order system with damping,
{\small \begin{equation}
\mathbf{\dot{x}}_i = A\mathbf{x}_i+ Bu_i, \quad \mathbf{x}_i(0)=\mathbf{x}_{0, i}  \quad \text{for} \quad i = 1, \hdots , N \label{n-agents}
 \end{equation}} 
\noindent We assume without loss of generality, 
{\small $A = \mbox{diag}(0,-b)$,
$B = [b^{-1},-b^{-1}]^{\top}$} and the damping coefficient $b>0$.
The system states are $\mathbf{x}_i= \begin{bmatrix} x_{1,i} &{x}_{2,i}\end{bmatrix}^{\top} \in \mathbb{R}^{2}$. The control input to agent $a_i$ satisfies the constraint $|u_i(t)|\le 1$ for all  $t\ge 0$. The fuel consumed by agent $a_i$ at time $t_f$ is $F(u_i(t),t_f)= \int_{0}^{t_f} |u_i(t)| dt$. In comparison to pure double integrator dynamics, this model includes the damping.
Let the maximum allowable fuel consumption be given as $\beta$; then the \emph{fuel budget constraint} is written as $F(u_i(t),t_f) \leq \beta$.
\begin{definition}
    A multi-agent system is said to achieve \emph{consensus} if for all $i,j \in \{1, \hdots, N\}$, where $i \neq j$, $||\mathbf{x}_i(t)-\mathbf{x}_j(t)|| \to 0$ as $t \to \overline{t}_f >0$ and $\mathbf{x}_i(t) = \mathbf{x}_j(t)$ for all $t \geq \overline{t}_f$. The time $\overline{t}_f$ is called the time to consensus, and the point $\mathbf{x}_i(\overline{t}_f)=\overline{\mathbf{x}} \mbox{ for }i = 1, \hdots, N$ is the corresponding consensus point. The consensus is said to be achieved in finite time if $\bar{t}_f<\infty$.
\end{definition} 

\begin{problem}
For a set of $N$ agents given by \eqref{n-agents},
compute the minimum time to consensus $\overline{t}_f$ and the corresponding consensus point $\mathbf{\overline{x}} \in \mathbb{R}^2$ such that $\mathbf{x}_i(\overline{t}_f)=\mathbf{x}_j(\overline{t}_f)=\mathbf{\overline{x}}$ and $F(u_i(t),\overline{t}_f) \leq \beta$ for all $i,j \in \{1,2 \hdots, N\}$ with $i \neq j$.
\end{problem}

\subsection{Preliminaries}
Let the set of admissible control inputs be defined as: 
\begin{equation}U_\beta:=\{u(t) \mbox{ s.t. } |u(t)|\le 1, F(u(t),t_f)\le \beta\}\label{admissible}\end{equation}
\begin{definition}
The \emph{Attainable Set} $\mathcal{A}_i^\beta(t_f,\mathbf{x}_{0,i})$ is the set of all the states that an agent $a_i$ can reach from $\mathbf{x}_{0,i} \in \mathbb{R}^2$ using admissible control $u_i(t) \in U_\beta$ at time $t_f>0$. 
Thus, {\small     $\mathcal{A}_i^\beta(t_f,\mathbf{x}_{0,i}):= \{ e^{At_f}\mathbf{x}_{0,i}+ \int_0^{t_f}\hspace{-0.1cm} e^{A(t_f-\tau)}Bu_i(\tau) d\tau,
        \forall u_i(t) \in U_\beta \}$ }
\end{definition}

Note that the set $\mathcal{A}^{\beta}(t_f,\mathbf{x}_{0})$ with $\mathcal{A}^{\beta}(t_f,\mathbf{0})$ are related by $\mathcal{A}^{\beta}(t_f,\mathbf{x}_{0}) = e^{At_f} \mathbf{x}_{0} + \mathcal{A}^{\beta}(t_f,\mathbf{0})$  \cite{fashoro1992reachability}.

\begin{definition}
\emph{Reachable Set} $\mathcal{R}_i^\beta(t_f,\mathbf{0})$ is the set of initial conditions from which an agent $a_i$ can reach the origin $\mathbf{0}\in \mathbb{R}^2$ using admissible control $u_i(t) \in U_{\beta}$ in time $t_f>0$. 
Thus, \small     $\mathcal{R}_i^\beta(t_f,\mathbf{0}):= \Big\{\int_0^{t_f}\hspace{-0.1cm} e^{-A\tau}Bu_i(\tau) d\tau,
        \forall u_i(t) \in U_{\beta} \Big\}$
    \normalsize
\end{definition}

  The $\mathcal{A}_i^\beta(t_f,\mathbf{x}_{0,i})$ is also obtained from $\mathcal{R}_i^\beta(t_f,\mathbf{0})$ by using $\mathcal{A}_i^\beta(t_f,\mathbf{x}_{0,i}) = e^{At_f}(\mathbf{x}_{i0}+\mathcal{R}_i^\beta(t_f,\mathbf{0}))$ \cite{fashoro1992reachability}. The points on the boundary of the set $\mathcal{R}_i^\beta(t_f,\mathbf{0})$ are those initial conditions which require $F(u_i(t),t_f)\ge \beta$ amount of fuel to reach the origin at time $t_f$. Thus, the control inputs that characterize the   $\mathcal{R}_i^\beta(t_f,\mathbf{0})$ are obtained by minimizing the functional $F(u_i(t),t_f)$. The following lemma from \cite{athans2013optimal,Hjek1979L1optimizationIL} gives the optimal control inputs.
\begin{lemma}\label{thm:bang-off-bang}  \cite{athans2013optimal,Hjek1979L1optimizationIL}
    For the system given by \eqref{n-agents}, control ${u}_i(t)$ that steers state-trajectory $\mathbf{x}_i(t)$ of agent $a_i$ from an initial condition $\mathbf{x}_{0,i} \in \mathbb{R}^2$ to the {  origin} with \emph{minimum} fuel i.e., $F(u_i(t),t_f)$ and in \emph{finite} time $t_f$ is of the form 
\begin{equation}
   {u}_i(t)=\begin{cases}
    \gamma_1,  & \text{$t \in [0,t_1]$}\\
    0, & \text{$t \in [t_1,t_2]$} \qquad  0 \leq t_1 \leq t_2 \leq t_f< \infty\\
    \gamma_2  & \text{$t \in [t_2,t_f]$}\\
  \end{cases}\label{twoswitchcontrol}
\end{equation}
where $\gamma_l \in \{+1,-1\}$ for $l=1,2$.  

\end{lemma}
Inputs \eqref{twoswitchcontrol} are represented as sequence $(\gamma_1,0,\gamma_2)$. The two fuel optimal control sequences from Lemma \ref{thm:bang-off-bang} are $s_1=(+1,0,-1)$ and $s_3=(-1,0,+1)$. These sequences satisfy the necessary conditions of Pontryagin's maximum principle \cite{pontryagin1987mathematical}. However, for \eqref{n-agents}, there are initial conditions for which the initial co-states are such that the optimal control exists, but is indeterminable. For these initial co-states the control sequences of the form $s_2=(-1,0,-1)$ and $s_4=(+1,0,+1)$ are also optimal \cite{Hjek1979L1optimizationIL}. These control sequences have a zero level making the control remain inactive for a finite duration, consuming no fuel. The fuel consumption for all the sequences up to time $t_f$ is then $F(u(t),t_f)=t_f-t_2+t_1.$ Since $\mathcal{A}_i^{\beta}(t_f,\mathbf{x}_{0,i}) = e^{At_f}(\mathbf{x}_{0,i}+\mathcal{R}_i^{\beta}(t_f,\mathbf{0}))$, the fuel optimal control sequences $s_p, p =1,2,3,4$ also characterize the attainable set  $\mathcal{A}_i^{\beta}(t_f,\mathbf{x}_{0,i})$. 
To characterize the boundary of the set $\mathcal{A}_i^{\beta}(t_f,\mathbf{x}_{0,i})$ we take {$F(u_i(t),t_f) = \beta$ which signifies} the full utilization of the fuel budget. Now, let
$\hat{U}_{\beta}:=\{u_i(t) \mbox{ s.t. \eqref{twoswitchcontrol} and } F(u_i(t),t_f)\le \beta\}\label{admissible}$. Then, the attainable set becomes $\mathcal{A}_i^{\beta}(t_f,\mathbf{x}_{0,i})
= \left\{e^{At_f}\mathbf{x}_{0,i}+\int_0^{t_f} e^{A(t_f-\tau)}Bu_i(\tau) d\tau, \forall u_i(t) \in \hat{U}_\beta \right\}$
    
\subsubsection{Attainable Set Characterization} \label{sec:attainableset}
The states obtained from $\mathbf{x}_{0,i}$ at time $t_f$ using the input profile \eqref{twoswitchcontrol} are 
{\small\begin{equation} 
    \mathbf{x}_i(t_f) = e^{At_f}\mathbf{x}_{0,i} + \Bigg( \gamma_1\int_0^{t_{1}} +~ \gamma_2\int_{t_{2}}^{t_{f}} \Bigg) e^{A(t_f-\tau)} B d\tau.
\label{eq:attainablestates}\end{equation}}Further, ensuring $t_f-t_2+t_1=\beta$ gives the boundary of the set $\mathcal{A}_i^\beta(t_f, \mathbf{x}_{0,i})$.
Expanding the equation \eqref{eq:attainablestates} for all sequences $s_1,s_2,s_3,s_4$, we obtain boundary expressions for $\mathcal{A}_i^\beta(t_f, \mathbf{x}_{0,i})$. For notation purpose, we denote $t_1,t_2$ as $t^{s_p}_1,t^{s_p}_2$ for sequences $s_p, 
\mbox{where~} p=1,2,3,4$ and drop the subscript $i$ for brevity. On expanding \eqref{eq:attainablestates}, we get both polynomial and exponential terms in expressions. To convert the exponential terms into polynomial terms, we substitute
\begin{equation}\label{eq:subs}q_f = e^{-bt_f}, w_1=e^{\frac{b^2}{2}x_1} , l_1 = e^{\frac{b}{2}(\beta - b x_{10})}, l_2 = e^{\frac{-b}{2}(\beta + b x_{10})}\end{equation}
Further, denote $L:=(l_1,l_2)$. Then for $s_1$ sequence the relation between $\mathbf{x}_{0} = [x_{10},{x}_{20}]$, $\mathbf{x}=[x_{1},{x}_{2}]$, $t_f$ and $\beta$ in terms of the substituted variables is $ \Gamma^{s_1}_{(L,x_{20})}(w_{1},x_{2},q_1)=0$ where {\small \begin{equation}
\label{eq:p10m1}
 \Gamma^{s_1}_{(L,x_{20})}(w_{1},x_{2},q_f)={x}_2-q_f{x}_{20} + \frac{1}{b^2} q_fw_1l_1 - \frac{q_f}{b^2} - \frac{1}{b^2}+ \frac{w_1l_2}{b^2}
\end{equation}} Also, the switching times in sequence $s_1$ are
$t^{s_1}_1 = \frac{ \beta + b(x_1 - x_{10})}{2}$ and $t^{s_1}_2 = t_f+ 
(\frac{ b(x_1 - x_{10})- \beta}{2})$. 
Similarly, for control input sequence $s_3$, we get $ \Gamma^{s_3}_{(L,x_{20})}(w_{1},x_{2},q_1)=0$
{\small  \begin{equation} \label{eq:m10p1}
\Gamma^{s_3}_{(L,x_{20})}(w_{1},x_{2},q_f)=w_1{x}_2 - w_1 q_f {x}_{20} - \frac{q_f}{b^2l_2}+\frac{q_f w_1}{b^2}+\frac{w_1}{b^2}-\frac{1}{l_1b^2} 
\end{equation}
}
and the switching times are $t^{s_3}_1 = \frac{ \beta - b(x_1 - x_{10})}{2}$,  $t^{s_3}_2 = t_f-( \frac{\beta + b(x_1 - x_{10})}{2})$.

For control input $s_2$, we directly have  $x_1 =  x_{10}-\frac{\beta}{b}$ and hence in terms of substituted variables $\Gamma^{s_2}_{(L,x_{20})}(w_{1},x_{2},q_f)=0$ with 
\begin{equation}\label{eqn:s2}
     \Gamma^{s_2}_{(L,x_{20})}(w_{1},x_{2},q_f) = w_1-\frac{1}{l_1}
\end{equation}
\begin{equation}
    x_2 =  q_fx_{20} + \frac{-q_f+1+q_f e^{bt_1}+ e^{-b(\beta - t_1)}}{b^2}
\end{equation}
and the switching time instances are
\begin{equation}
    t_1^{s_2} \;=\; \frac{1}{b} \, 
\ln \!\left( 
\frac{x_2b^2-q_fx_{20}b^2+q_f-1}{q_f - e^{-b \beta}} 
\right)
\end{equation}
\begin{equation}
t_2^{s_2} \;=\; t_f - \beta 
\;+ t_1^{s_2}
\end{equation}
Similarly for the sequence $s_4$, we get $x_1 =  x_{10}+\frac{\beta}{b}$ and thus, $\Gamma^{s_4}_{(l,x_{20})}(w_{1},x_{2},q_f)=0$
\begin{equation}\label{eqn:s4}
    \Gamma^{s_4}_{(l,x_{20})}(w_{1},x_{2},q_f) = w_1-\frac{1}{l_2}
\end{equation}
\begin{equation}
    x_2 = q_fx_{20} + \frac{e^{bt_1}}{w_{10}^2l_1^2b^2} + \frac{q_f-1-q_f e^{bt_1}}{b^2} 
\end{equation}
where, $w_{10}=e^{-\frac{b^2}{2}x_{10}}$
and \begin{equation}
t_1^{s_4} \;=\; \frac{1}{b} \, 
\ln \!\left( 
\frac{x_2b^2-q_fx_{20}b^2-q_f+1}{-q_f + e^{-b\beta}} 
\right)
\end{equation}
\begin{equation}
t_2^{s_4} \;=\; t_f - \beta +t_1^{s_4}.
\end{equation}
The set of points $\mathbf{x}\in\mathbb{R}^2$ that satisfy equation \eqref{eq:p10m1} or \eqref{eq:m10p1} forms a part of the boundary $\partial\mathcal{A}^\beta(t_f,\mathbf{x}_{0})$ for the two control sequences $s_1$ and $s_3$ respectively. Let us denote these two parts by the following sets
{\small\begin{equation}  \label{eq:bds1}
\partial{\mathcal{A}}^{\beta,s_1}=\{ \mathbf{x}~\mbox{s.t.} ~ \Gamma^{s_1}_{(l,x_{20})}=0  \text{, } 0\leq t^{s_1}_1 \leq t^{s_1}_2 \le t_f < \infty \}
\end{equation}
\begin{equation} \label{eq:bds3} 
\partial{\mathcal{A}}^{\beta,s_3}:= \{ \mathbf{x}~\mbox{s.t.} ~ \Gamma^{s_3}_{(l,x_{20})}=0 \text{, } 0 \leq  t^{s_3}_1 \leq  t^{s_3}_2 \le t_f < \infty 
\}
\end{equation}}

From equations \eqref{eqn:s2} and \eqref{eqn:s4} respectively, we see that this part of the boundary has a fixed ${x}_1$ coordinate. Thus, 
{\small \begin{equation}\label{eq:bds2} \partial{\mathcal{A}}^{{\beta},s_2}=\bigl\{ \mathbf{x}~\mbox{s.t.} ~ \Gamma^{s_2}_{(l,x_{20})}=0  ,0 \le t_1^{s_2} \leq t_2^{s_2} \le t_f < \infty\bigr\}
\end{equation}
 \begin{equation} \label{eq:bds4} \partial{\mathcal{A}}^{{\beta},s_4}=\bigl\{ \mathbf{x}~\mbox{s.t.} ~ \Gamma^{s_4}_{(l,x_{20})}=0  ,0 \le t_1^{s_4} \leq t_2^{s_4}\le t_f < \infty\bigr\}
 \end{equation}}

\begin{figure}[h]
 \centering 
{\includegraphics[width=.7\linewidth, height=3.8cm]{{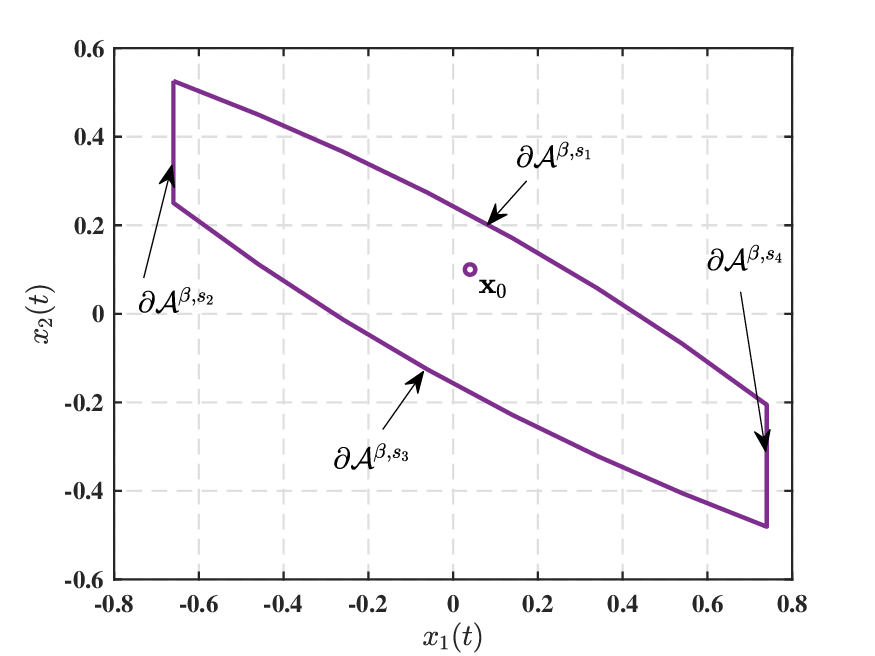}}%
 \caption{Attainable set }\label{fig:a}}
\end{figure}

\noindent The boundary of the set is $\mathcal{A}^{\beta}(t_f,\mathbf{x}_{0})=\bigcup_{p=1}^4\partial{\mathcal{A}}^{{\beta},s_p}$. An example of attainable set with $t_f = 1.4918$ and $\beta=0.7$ is shown in Figure \ref{fig:a}.

We list several important properties that hold true for $\mathcal{A}^{\beta}(t_f,\mathbf{x}_{0})$, which are already explored in \cite{tacsubmitted2024} in the form of Lemma \ref{lem:propertiesofattainableset} (see also c.f. \cite{Hjek1979L1optimizationIL} for similar properties for reachable set).

\begin{lemma} (c.f. \cite{tacsubmitted2024,sclsubmitted2025})\label{lem:propertiesofattainableset}
For the attainable set $\mathcal{A}^{\beta}(t_f, \mathbf{x}_{0})$, 
the following structural properties holds true:

(i) For every $t_f>0$ and $\beta \in (0,\infty]$, the set 
    $\mathcal{A}^{\beta}(t_f,\mathbf{x}_{0})$ is a convex subset of $\mathbb{R}^2$.

(ii)  If the fuel budget satisfies $\beta \geq t_f$, then
\[
\mathcal{A}^{\beta}(t_f,\mathbf{x}_{0}) \;=\; \mathcal{A}^\infty(t_f,\mathbf{x}_{0}).
\]

(iii)   For any finite fuel budgets $\beta', \beta \in (0,\infty)$ with $\beta' \leq \beta$,
    \[
    \mathcal{A}^{\beta'}(t_f,\mathbf{x}_{0})
    \;\subseteq\;
    \mathcal{A}^{\beta}(t_f,\mathbf{x}_{0})
    \;\subseteq\;
    \mathcal{A}^\infty(t_f,\mathbf{x}_{0}).
    \]
\end{lemma}

\begin{lemma}
As the time $t_f \to \infty$, the attainable set converges 
    in the set-wise sense to $\mathcal{A}^{\beta}(\infty, \mathbf{x}_{0})$, i.e.,
    \[
    \lim_{t_f \to \infty} \mathcal{A}^{\beta}(t_f,\mathbf{x}_{0}) 
    \;=\; \mathcal{A}^{\beta}(\infty, \mathbf{x}_{0}).
    \] where, $\mathcal{A}^{\beta}(\infty, \mathbf{x}_{0}) = \mbox{conv}(\partial\mathcal{A}^{\beta,s_1}(\infty, \mathbf{x}_{0}), \partial{\mathcal{A}^{\beta,s_3}(\infty, \mathbf{x}_{0})})$, $\partial\mathcal{A}^{\beta,s_1}(\infty, \mathbf{x}_{0}):= \{ (x_1, x_2) \mbox{ s.t. } x_1 \in [ {x}_{10} - \frac{\beta}{b},  {x}_{10} + \frac{\beta}{b}], x_2 = \frac{1}{b^2}(1 - e^{\frac{b^2}{2}x_1}l_2) \}$, and $\partial\mathcal{A}^{\beta,s_3}(\infty, \mathbf{x}_{0}):= \{ (x_1, x_2) \mbox{ s.t. } x_1 \in [ {x}_{10} - \frac{\beta}{b},  {x}_{10} + \frac{\beta}{b}], x_2 = \frac{1}{l_1 b^2}(e^{\frac{-b^2}{2}x_1}-l_1) \}$. Here $\mbox{conv}$ is convex hull.
    
\end{lemma}

\begin{proof}
    The lemma follows by substituting $q_f = 0$ in \eqref{eq:p10m1} and \eqref{eq:m10p1}.
\end{proof}

\section{\label{subsec:roc} Region of Consensus}
\begin{figure}[h]
 \centering 
 {\includegraphics[width=.7\linewidth, height=4cm]{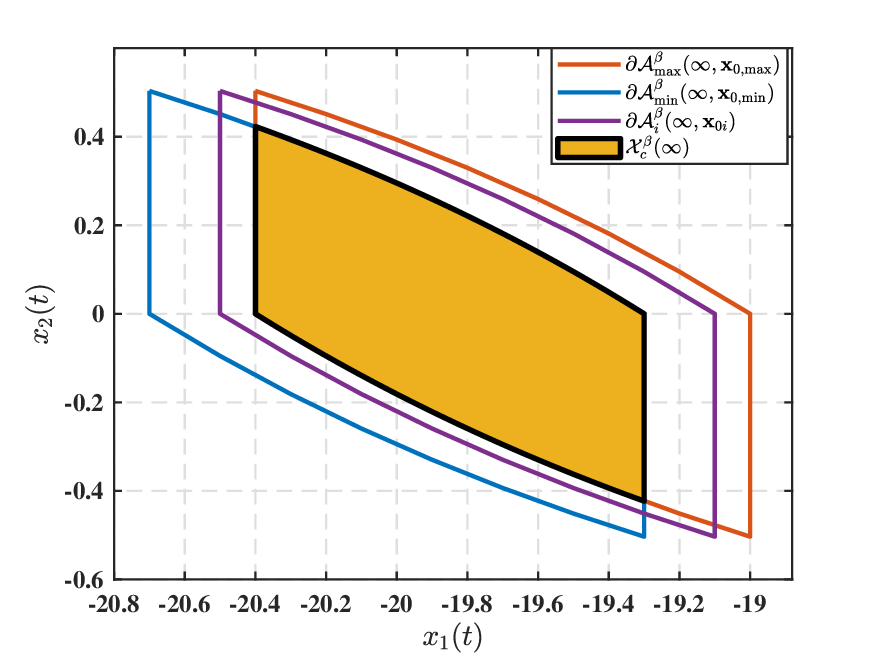}%
 \caption{Region of Consensus }\label{fig:rcb}}
\end{figure}

\noindent For consensus to be possible, we need{\small $\bigcap\limits_{i=1}^N{\mathcal{A}}_{i}^\beta(\tau, \mathbf{x}_{0,i})\ne\emptyset$}. Now, let ${x}_{1,\max}:= \max\limits_{i\in\{1,..,N\}}\,{x}_{10,i}$ and ${x}_{1,\min}:= \min\limits_{i\in\{1,..,N\}}\,{x}_{10,i}$ denote, respectively, the maximum and minimum ${x}_1$ co-ordinates among the initial conditions of all agents $a_i,~ i=1,..,N$. Also, let the corresponding initial conditions be labeled as $\mathbf{x}_{0,\max}$ and $\mathbf{x}_{0,\min}$ respectively. Accordingly, denote the corresponding attainable sets as $\mathcal{A}_{\max}^{\beta}(t_f,\mathbf{x}_{0,\max})$ and $\mathcal{A}_{\min}^{\beta}(t_f,\mathbf{x}_{0,\min})$ respectively.
 \noindent  Let $\mathcal{X}^{\beta}_c(t_f) := \mathcal{A}_{\max}^{\beta}(t_f,\mathbf{x}_{0,\max}) \cap \mathcal{A}_{\min}^{\beta}(t_f,\mathbf{x}_{0,\min})$. For some $t_f<\bar{t}_f$, $\mathcal{X}^{\beta}_c(t_f)$ is empty. Also, as $t_f\to \infty$, the set converges to $\mathcal{X}_c^{\beta}(\infty)$. The boundary of set $\mathcal{X}_c^{\beta}(\infty)$ consists of two portions $\partial\mathcal{X}_c^{\beta,s_1}(\infty, \mathbf{x}_{0,\min}):= \{ (x_1, x_2) \mbox{ s.t.~} x_1 \in [ {x}_{10,\max} - \frac{\beta}{b},  {x}_{10,\min} + \frac{\beta}{b}], x_2 = \frac{1}{b^2}(1 - e^{\frac{b^2}{2}x_1}l_{2,\min}) \}$ and $\partial\mathcal{X}_c^{\beta,s_3}(\infty, \mathbf{x}_{0,\max}):= \{ (x_1, x_2) \mbox{ s.t.~} x_1 \in [ {x}_{10,\max} - \frac{\beta}{b},  {x}_{10,\min} + \frac{\beta}{b}],  x_2 = \frac{1}{l_{1,\max} b^2}(e^{\frac{-b^2}{2}x_1}-l_{1,\max}) \}$. Then, the entire set $\mathcal{X}_c^{\beta}(\infty)$ is the following convex hull (for e.g. see Fig. \ref{fig:rcb}). {\small\begin{equation}\label{consensuspoints}
\mathcal{X}_c^{\beta}(\infty)=\mbox{conv}(\partial\mathcal{X}_c^{\beta,s_1}(\infty,\mathbf{x}_{0,\min}), \partial{\mathcal{X}_c^{\beta,s_3}(\infty,\mathbf{x}_{0,\max})})
\end{equation}} 
\begin{theorem} \label{thm:consensuspossible}
Consensus is possible if and only if the set $\mathcal{X}_c^{\beta}(\infty)$ is non-empty. The consensus point $\bar{\mathbf{x}}$ always satisfies $\bar{\mathbf{x}}\in \mathcal{X}_c^{\beta}(\infty)$.
\end{theorem}
\begin{proof}
It is clear that if $\mathcal{X}_c^{\beta}(\infty)=\emptyset$ then consensus is not possible. Further if $\mathcal{X}_c^{\beta}(\infty) \ne \emptyset$, then since the set $\mathcal{X}_c^{\beta}(\infty)$ is contained in the attainable set $\mathcal{A}_{i}^{\beta}(\infty,\mathbf{x}_{0,i})$ for all agents $i=1,..,N$. Thus, consensus is possible. 
\end{proof}

This gives us the following corollary.
\begin{corollary}
Consensus is possible if and only if the initial conditions of all agents are such that 
\begin{equation} \label{eq:consensuspossible}
{x}_{1,\max}- {x}_{1,\min} \le \frac{2\beta}{b} \end{equation}
\end{corollary}
\begin{proof}
The proof follows from Theorem \ref{thm:consensuspossible} and the fact that the set $\mathcal{X}_c^{\beta}(\infty)$ is non-empty if and only if ${x}_{1,\max}-\frac{\beta}{b} \leq {x}_{1,\min}+\frac{\beta}{b}$.  
\end{proof}

\begin{remark}
    Note that if ${x}_{1,\max}- {x}_{1,\min} < \frac{2\beta}{b}$ then the set $\mathcal{X}_c^{\beta}(\infty)$ is of non-zero Lebesgue measure. Thus, the set $\mathcal{X}_c^{\beta}(t_f)$ also has non-zero measure at some finite time $t_f<\infty$. Thus,  if ${x}_{1,\max}- {x}_{1,\min} < \frac{2\beta}{b}$ then, minimum time to consensus will also be finite i.e., $\bar{t}_f<\infty$. On the other hand, if ${x}_{1,\max}- {x}_{1,\min}=\frac{2\beta}{b}$,  then $\mathcal{X}^{\beta}_c(\infty) = {({x}_{1,\max}-\frac{\beta}{b},0)}$ and for any finite time $t_f<\infty$, $\mathcal{X}_c^{\beta}(t_f)$ is empty. Therefore, consensus in finite time is not possible.   
\end{remark}

\begin{assumption}\label{assumption1:cp}
    We assume that consensus is possible for some $t_f<\infty$ and thus, the initial conditions of all agents satisfy ${x}_{1,\max}- {x}_{1,\min}< \frac{2\beta}{b}$.
\end{assumption}
\section{\label{sec:min.timeconsensus}Minimum Time Consensus}
Our goal is to find {\small $\bar{t}_f:=\min\limits_{\bigcap\limits_{\alpha=1}^N{\mathcal{A}}_{\alpha}^\beta(\tau, \mathbf{x}_{0,\alpha})\ne\emptyset}{\tau}$} and the corresponding {\small $\bar{\mathbf{x}}\in\bigcap\limits_{\alpha=1}^N{\mathcal{A}}_{\alpha}^\beta(\bar{t}_f, \mathbf{x}_{0,\alpha})$}. For this we require $\bigcap\limits_{\alpha=1}^N{\mathcal{A}}_{\alpha}^\beta(\bar{t}_f, \mathbf{x}_{0,\alpha})$ which is combinatorially intractable. However, note that  
$\mathcal{A}_i^\beta(t_f, \mathbf{x}_{0,i})$ is a convex set. Thus, to make the computation of the minimum time consensus point tractable, we make use of Helly's theorem, stated next.
\begin{theorem}{\emph{Helly's Theorem}~\cite{gruber2007convex}}:
    Let $F$ be a finite family of at least $n+1$ convex sets on $\mathbb{R}^n$. If the intersection set of every $n+1$ members of $F$ is non-empty, then all members of $F$ have a non-empty intersection.
\end{theorem}

In our problem, $n=2$, the computation of the minimum-time consensus point can be reduced to computing the minimum-time consensus for all possible triplets of agents.
\subsection{Minimum time consensus for a triplet}

 {Let the minimum time to consensus for a triplet $\{a_i,a_j,a_k\}$ be {\small $\bar{t}_{ijk} := \min\limits_{\bigcap\limits_{\alpha\in\{i,j,k\}}{\mathcal{A}}_{\alpha}^\beta(\tau, \mathbf{x}_{0,\alpha})\ne\emptyset}{\tau}$}. Since, we need {\small $\bigcap\limits_{\alpha\in\{i,j,k\}}{\mathcal{A}}_{\alpha}^\beta(\tau, \mathbf{x}_{0,\alpha})\ne\emptyset$}, for every pair of agents, $(a_i,a_j)$,$(a_j,a_k)$ and $(a_i,a_k)$, their attainable sets must intersect first. Let $\bar{t}_{ij}:=\min\limits_{\bigcap\limits_{\alpha\in\{i,j\}}{\mathcal{A}}_{\alpha}^\beta(\tau, \mathbf{x}_{0,\alpha})\ne\emptyset}{\tau}$. Similarly, $\bar{t}_{jk}$ and $\bar{t}_{ik}$ are defined for $(a_j,a_k)$ and $(a_i,a_k)$. Then, $\bar{t}_{ijk} \geq \text{max}\{\bar{t}_{ij} ,\bar{t}_{jk},\bar{t}_{ik}\}$.
 WLOG $\bar{t}_{ij} = \text{max}\{\bar{t}_{ij} ,\bar{t}_{jk},\bar{t}_{ik}\}$ and the corresponding point of intersection $\bar{\mathbf{x}}_{ij}\in{\mathcal{A}}_{i}^\beta(\bar{t}_{ij}, \mathbf{x}_{0,i})\cap\mathcal{A}_{j}^\beta(\bar{t}_{ij}, \mathbf{x}_{0,j})$.} 
We now have two cases: 1) $\bar{\mathbf{x}}_{ij} \in \mathcal{A}_k^\beta(\bar{t}_{ij},\mathbf{x}_{0,k})$
and  2) $\bar{\mathbf{x}}_{ij} \notin \mathcal{A}_k^\beta(\bar{t}_{ij},\mathbf{x}_{0,k})$

\subsubsection*{Case 1} We need to compute the solution to the expressions defining the boundary of two attainable sets and check if it is in the third agent's attainable set. This is done by checking if the point $\bar{\mathbf{x}}_{ij}$ satisfies inequalities associated with the boundary of $\mathcal{A}_k^\beta(\bar{t}_{ij},\mathbf{x}_{0,k})$.
Computing $\bar{t}_{ij}:={\min\limits_{\substack{\bigcap\limits_{\alpha\in\{i,j\}}{\mathcal{A}}_{\alpha}^\beta(\tau, \mathbf{x}_{0,\alpha})\ne\emptyset}} \tau}$ along with $\bar{\mathbf{x}}_{ij}$ requires the intersection $\bigcap\limits_{\nu\in\{i,j\}}{\mathcal{A}}_{\nu}^\beta(\tau, \mathbf{x}_{0,\nu})\ne\emptyset$. This is possible in $2$ ways i.e., $\partial\mathcal{A}^{\beta,s_1}_i\cap \partial\mathcal{A}^{\beta,s_3}_j$ or $\partial\mathcal{A}^{\beta,s_3}_i\cap \partial\mathcal{A}^{\beta,s_1}_j$. Now  $\partial\mathcal{A}^{\beta,s_p}_i$ for $p=1,2,3,4$ are given in \eqref{eq:bds1}, \eqref{eq:bds2}, \eqref{eq:bds3} and \eqref{eq:bds4}. Further, due to substitution \eqref{eq:subs}, the boundary is in terms of $q_1 = e^{-b\tau}$, $w_1 = e^{\frac{b^2}{2}x_1}$ and $x_2$. By maximizing $q_1$ subject to the constraints of $\partial\mathcal{A}^{\beta,s_1}_i\cap \partial\mathcal{A}^{\beta,s_3}_j$ or $\partial\mathcal{A}^{\beta,s_3}_i\cap \partial\mathcal{A}^{\beta,s_1}_j$  gives us $\bar{t}_{ij}$. This maximization is done by solving simultaneously $\Gamma^{s_1}_{(L_i,x_{20,i})}({w}_{1},{x}_{2},q_{1}) = \Gamma^{s_3}_{(L_j,x_{20,j})}({w}_{1},{x}_{2},q_{1}) = 0$ and  $\Gamma^{s_3}_{(L_i,x_{20,i})}({w}_{1},{x}_{2},q_{1}) = \Gamma^{s_1}_{(L_j,x_{20,j})}({w}_{1},{x}_{2},q_{1}) = 0$. In both situations, the solution is obtained by using Gr\"{o}bner basis with lexicographic ordering of $({x}_{2}\succ {w}_{1}  \succ q_{1})$. This allows us to express $q_1$ as a rational function of $w_1$ which we maximize w.r.t $w_1$. Finally, the largest of the two solutions gives $\bar{t}_{ij}$ and the corresponding intersection point $\bar{\mathbf{x}}_{ij}$. If  $\bar{\mathbf{x}}_{ij} \in \mathcal{A}_k^\beta(\bar{t}_{ij},\mathbf{x}_{0,k})$ then $\bar{t}_{ij}=\bar{t}_{ijk}$ and $\bar{\mathbf{x}}_{ij}= \bar{\mathbf{x}}_{ijk}$. However, if the condition of case 1 is invalid, then we proceed to case 2.


 \subsubsection*{Case 2} In case 2, even though for each pair of agents, the corresponding attainable sets intersect at time $\bar{t}_{ij} = \text{max}\{\bar{t}_{ij} ,\bar{t}_{jk},\bar{t}_{ik}\}$, but the intersection of all the three attainable sets is empty at $\bar{t}_{ij}$. 
Therefore, the first point of intersection will always occur on the boundaries of the attainable sets all the three agents.
We note that an intersection must exist at some time $\bar{t}_{ijk}>\bar{t}_{ij}$ on the boundaries of all three attainable sets. Therefore, we must obtain all  solutions $t_f$ that satisfy  
 \begin{equation}\label{intersectionbd}
\partial {\mathcal{A}}_i^\beta(t_f, \mathbf{x}_{0,i}) \cap \partial {\mathcal{A}}_j^\beta(t_f, \mathbf{x}_{0,j}) \cap \partial{\mathcal{A}}_k^\beta(t_f, \mathbf{x}_{0,k}) \neq \emptyset.
\end{equation}
Then $\bar{t}_{ijk}$ is chosen to be the least $t_f$ among all solutions. Further, $\mathbf{\bar{x}}_{ijk} \in \partial {\mathcal{A}}_i^\beta(\bar{t}_{ijk}, \mathbf{x}_{0,i}) \cap \partial{\mathcal{A}}_j^\beta(\bar{t}_{ijk}, \mathbf{x}_{0,j})\cap \partial{\mathcal{A}}_k^\beta(\bar{t}_{ijk}, \mathbf{x}_{0,k})$ is the only point in the intersection.
{  Since $\mathbf{\bar{x}}_{ijk}$ lies on the boundary of attainable sets, control inputs of the form \eqref{twoswitchcontrol} are required to drive the states of all agents to $\mathbf{\bar{x}}_{ijk}$ at time $\bar{t}_{ijk}$.  Therefore, depending upon the structure of the input, $\mathbf{\bar{x}}_{ijk}$ may lie on one of the four boundaries of $\partial{\mathcal{A}}_i^\beta(t_f, \mathbf{x}_{0,i})$, $\partial{\mathcal{A}}_j^\beta(t_f, \mathbf{x}_{0,j})$ and $\partial{\mathcal{A}}_k^\beta(t_f, \mathbf{x}_{0,k})$.  Solving equation \eqref{intersectionbd} requires computing solutions of a set of polynomial equations with inequality constraints simultaneously. Since there are $4$ parts to the boundaries, there are $64$ sets of three equations to be solved for computing $\mathbf{\bar{x}}_{ijk}$. Using arguments and using Assumption 1 discussed in \cite[Section III-A]{tacsubmitted2024}, we only need to consider $18$ cases to compute the minimum time consensus point. 
Among all the 18 possibilities, we discuss the solution of these polynomial equations for specific possibilities.

\subsubsection{$\partial{\mathcal{A}}_{i}^{\beta,s_1} \cap \partial{\mathcal{A}}_{j}^{\beta,s_3} \cap \partial{\mathcal{A}}_{k}^{\beta,s_1} $  }\label{sec:boundaryequations1}

 In the boundary equations for $\partial {\mathcal{A}}_{i}^{\beta,s_1}$, $\partial {\mathcal{A}}_{j}^{\beta,s_3}$ and $\partial {\mathcal{A}}_{k}^{\beta,s_1}$, we substitute $\mathbf{x}_i = \mathbf{x}_j = \mathbf{x}_k = \hat{\mathbf{x}}_{1}$. We are solving the boundary equations for $t_{f,1}$, and $\hat{\mathbf{x}}_1$ which has two components i.e., $\hat{\mathbf{x}}_1=[\hat{x}_{1,1}, \hat{x}_{2,1}]^{\top}$. Also, recall from the substitution \eqref{eq:subs} $w_{1,1}=e^{\frac{b^2}{2}x_{1,1}}$ and $q_{f,1}=e^{-bt_{f,1}}$.   As a result, we solve the three equations $\Gamma^{s_1}_{(L_i,x_{20,i})}(\hat{w}_{1,1},\hat{x}_{2,1},q_{f,1}) = \Gamma^{s_3}_{(L_j,x_{20,j})}(\hat{w}_{1,1},\hat{x}_{2,1},q_{f,1}) = \Gamma^{s_1}_{(L_k,x_{20,k})}(\hat{w}_{1,1},\hat{x}_{2,1},q_{f,1})=0$ simultaneously for $q_{f,1}$, $\hat{{w}}_{1,1}$ and $\hat{{x}}_{2,1}$ using the Gr\"{o}bner Basis based elimination \cite{cox2013ideals}. Since there are three equations and three unknowns $q_{f,1}$, $\hat{{w}}_{1,1}$ and $\hat{{x}}_{2,1}$, we use lexicographic ordering $(\hat{{w}}_{1,1}\succ \hat{{x}}_{2,1} \succ q_{f,1})$ to arrive at the final expressions. The computation of the Gr\"{o}bner basis is performed using the Sagemath \cite{stein2007sage}.  
Subsequently, we get a cubic equation in 
$b_3q_{f,1}^3+b_2q_{f,1}^2+b_1q_{f,1}+b_0=0$.
Among all the real solutions to the cubic equation, we choose the one with the lowest real value in the interval $(0,1)$ as  $q_{f,1}$. Then, the expression of time is $t_{f,1} =  -\log(q_{f,1})$. The corresponding expressions of states are given as:

{\footnotesize \begin{align}\label{eq:st01}
\hat{w}_{1,1} &= \frac{-(\hat{x}_{2,1} - q_{f,1} x_{20,k} - q_{f,1} - 1))}{q_{f,1} l_{1,k} + l_{2,k}}\nonumber \\
\hat{x}_{1,1}&  = \frac{2}{b^2}\log(\hat{w}_{1,1})\\ 
\hat{x}_{2,1} & = -\frac{\left(\splitfrac{q_{f,1}^2(x_{20,i}l_{1,k} -l_{1,i}x_{20,k} - l_{1,i} + l_{1,k}) - l_{2,i} + l_{2,k}}{ +q_{f,1}(x_{20,i}l_{2,k} - l_{1,i} - l_{2,i}x_{20,k} - l_{2,i} + l_{1,k} + l_{2,k}) }\right)}{ q_{f,1}(l_{1,i} - l_{1,k}) + l_{2,i} - l_{2,k} }
\end{align}}
At the computed time $t_{f,1}$, these sets meet at a unique point $\hat{\mathbf{x}}_1$. An example of a unique intersection point with $\partial{\mathcal{A}}_{i}^{\beta,s_1} \cap \partial{\mathcal{A}}_{j}^{\beta,s_3} \cap \partial{\mathcal{A}}_{k}^{\beta,s_1} $ is shown in Fig.\ref{fig2}.
\begin{figure}[h!]
\centering 
{\includegraphics[width=.7\linewidth, height=4cm]{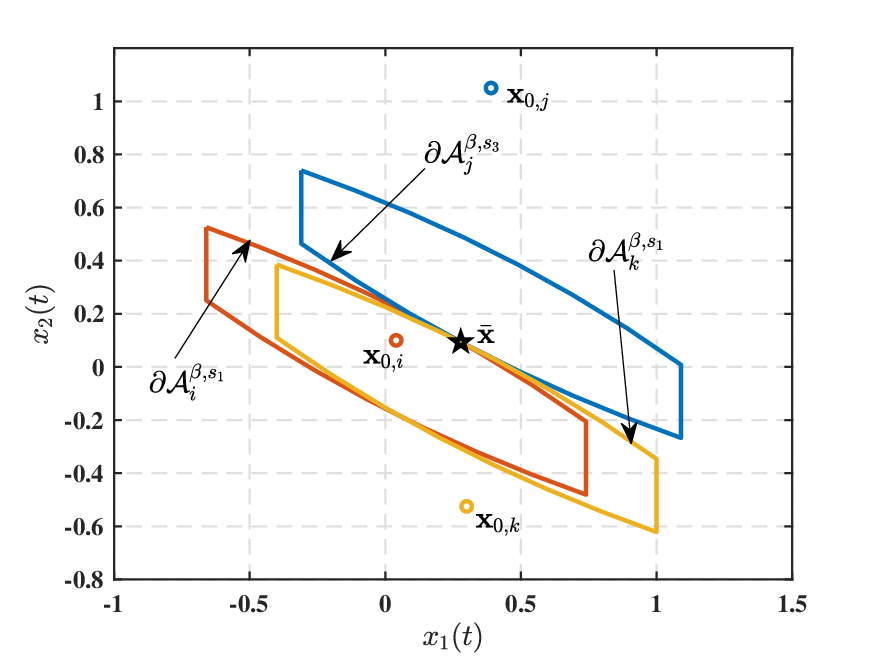}}
\caption{Intersection of Attainable sets}  \label{fig2}
\end{figure}

\subsubsection{Scenario 7: $\partial{\mathcal{A}}_{i}^{\beta_2,s_3} \cap \partial{\mathcal{A}}_{j}^{\beta_2,s_2} \cap \partial{\mathcal{A}}_{k}^{\beta_2,s_1}$} 
Similarly, we simultaneously solve equations, $\Gamma^{s_3}_{(L_i,x_{20i})}(\hat{w}_{1,7},\hat{x}_{2,7},q_{f,7}) = \Gamma^{s_2}_{(L_j,x_{20j})}(\hat{w}_{1,7},\hat{x}_{2,7},q_{f,7}) = \Gamma^{s_1}_{(L_k,x_{20k})}(\hat{w}_{1,7},\hat{x}_{2,7},q_{f,7})=0$
for $q_{f,7}$, $\hat{{w}}_{1,7}$ and $\hat{{x}}_{2,7}$. 

The equation $\Gamma^{s_2}_{(L_j,x_{20j})}(\hat{w}_{1,7},\hat{x}_{2,7},q_{f,7})=0$ gives us $w_{1,7}=\frac{1}{l_{1j}}$. Then, from the substitution \eqref{eq:subs}, we get $\hat{x}_{1,7} = x_{10j}-\frac{\beta}{b}$.
The expression for $t_{f,7}$ and the corresponding states obtained from

{ \small \begin{equation}
\begin{aligned}\label{eq:tf07}
q_{f,7} & = \frac{-(\hat{w}_{1,7}^2l_{2k}  - 2\hat{w}_{1,7} + l_{1i})}{\hat{w}_{1,7}^2l_{1k} + \hat{w}_{1,7}x_{20i} - \hat{w}_{1,7}x_{20k} - 2\hat{w}_{1,7} + l_{2i} } \\
t_{f,7} &= -\log(q_{f,7})
\end{aligned}
\end{equation}
}
{\small \begin{equation}
\begin{aligned} \label{eq:st07}
\hat{x}_{2,7} & = - \hat{w}_{1,7} q_{f,7}l_{1k} - \hat{w}_{1,7} l_{2k} + q_{f,7}x_{20k} + q_{f,7} + 1\\
 \end{aligned}
 \end{equation}
} 
The expressions of $t_{f,m}$ and $\hat{\mathbf{x}}_m$ for the remaining possibilities (2 to 6 and 8 to 18) can be obtained similarly. 

Now, the control input that can drive agent $\alpha\in\{i,j,k\}$ from its initial state to $\hat{\mathbf{x}}_m$ at time $t_{f,m}$ is of the form \eqref{twoswitchcontrol}. The expressions for the switching time instances, i.e. $t^{m,s_1}_{1\alpha}$, $t^{m,s_1}_{2\alpha}$ and $t^{m,s_3}_{1\alpha}$, $t^{m,s_3}_{2\alpha}$ for each agent $\alpha$ are obtained as follows. For $m=1,2,...,18$, if $\hat{\mathbf{x}}_m\in\partial{\mathcal{A}}_{\alpha}^{\beta,s_1}$, then the switching time instances for $a_\alpha$ are given as:

{\small \begin{equation}
 \begin{aligned} \label{eq:sw1}
 t_{1\alpha}^{m,s_1} & = \frac{ \beta + b(\hat{x}_{1,m} - x_{10 \alpha})}{2} \\
 t_{2\alpha}^{m,s_1} & = t_f+ 
\frac{ b(\hat{x}_{1,m} - x_{10\alpha})- \beta}{2}.
\end{aligned}\end{equation}}
\noindent Likewise, if $\hat{\mathbf{x}}_m\in\partial{\mathcal{A}}_{\alpha}^{\beta,s_3}$ then the switching time instances are
{\small \begin{equation} \begin{aligned}\label{eq:sw3}
    t_{1 \alpha}^{m,s_3} & = \frac{ \beta - b(\hat{x}_{1,m} - x_{10\alpha})}{2} \\ 
     t_{2 \alpha}^{m,s_3} & = t_f- \frac{(\beta + b(\hat{x}_{1,m} - x_{10\alpha}))}{2}.
\end{aligned} \end{equation}
} For $\hat{\mathbf{x}}_m \in \partial{\mathcal{A}}_{\alpha}^{\beta,s_2}$, 
\begin{equation}
\begin{aligned}\label{eq:sw2}
    t_{1 \alpha}^{m,s_2} & = \frac{1}{b} \, 
\ln \!\left( 
\frac{\hat{x}_{2,m}b^2-q_{f,m}x_{20 \alpha}b^2+q_{f,m}-1}{q_{f,m} - e^{-b \beta}} 
\right)\\ 
     t_{2 \alpha}^{m,s_2} & =  t_{f,m} - \beta 
\;+ t_{1 \alpha}^{m,s_2}.
\end{aligned}
\end{equation}
Lastly if $\hat{\mathbf{x}}_m \in \partial{\mathcal{A}}_{\alpha}^{\beta,s_4}$,
\begin{equation}
\begin{aligned}\label{eq:sw4}
    t_{1 \alpha}^{m,s_4} & = \frac{1}{b} \, 
\ln \!\left( 
\frac{\hat{x}_{2,m}b^2-q_{f,m}x_{20,\alpha}b^2-q_{f,m}+1}{-q_{f,m} + e^{-b\beta}} 
\right)\\ 
     t_{2 \alpha}^{m,s_4} & =  t_{f,m} - \beta 
\;+ t_{1 \alpha}^{m,s_4}.
\end{aligned}
\end{equation}
Among all the possibilities, we pick the value of $m$, denoted by $\bar{m}$, for which the inequality $0 \leq t^{m,\bullet}_{1 \alpha} \leq t^{m,\bullet}_{2 \alpha} \leq t_{f,m} < \infty$ holds for all $\alpha \in \{ i,j,k\}$, and $t_{f,m}$ is the \emph{least}. Then, the \emph{minimum time} to consensus is $\bar{t}_{ijk} = t_{f,\bar{m}}$ and the corresponding value of $\bar{\mathbf{x}}_{ijk}=\hat{\mathbf{x}}_{\bar{m}}$ is the minimum time consensus point for a triplet of agents $a_i,a_j,a_k$. 

\subsection{Minimum time consensus for N agents}
Now, for $N$-agent case, we need the following results which follow using exactly the same arguments as given in \cite{tacsubmitted2024,sclsubmitted2025}.
\begin{lemma} c.f. \cite{tacsubmitted2024,sclsubmitted2025}\label{intersection}
Consider a triplet of agents $\{a_i,a_j,a_k\}$ with initial conditions satisfying the inequality \eqref{eq:consensuspossible} and the time ${t}_f>0$ s.t. $\mathcal{A}_i^\beta({t}_f, \mathbf{x}_{0,i}) \cap \mathcal{A}_j^\beta({t}_f, \mathbf{x}_{0,j}) \cap \mathcal{A}_k^\beta({t}_f, \mathbf{x}_{0,k}) \neq \emptyset$ then for all time $t_f'>{t}_f$ the intersection of their attainable sets remains non empty i.e., $\mathcal{A}_i^\beta(t_f', \mathbf{x}_{0,i}) \cap \mathcal{A}_j^\beta(t_f', \mathbf{x}_{0,j}) \cap \mathcal{A}_k^\beta(t_f', \mathbf{x}_{0,k}) \neq \emptyset.$  
\end{lemma}


From Helly's theorem and Lemma \ref{intersection}, we get the following result. 

\begin{theorem} c.f. \cite{tacsubmitted2024,sclsubmitted2025}\label{thm:n-consensus}
     Consider $N$ agents with initial conditions satisfying the inequality \eqref{eq:consensuspossible}. Then the minimum time to consensus $\bar{t}_f = \max\limits_{1 \leq i,j,k \leq N}\bar{t}_{ijk}.$ 
\end{theorem}

 Theorem \ref{thm:n-consensus} gives us the minimum time to consensus $\bar{t}_f$ and the corresponding consensus point $\bar{\mathbf{x}}$ using all possible $\Mycomb[N]{3}$ triplet of agents. Assuming all agents are communicating over a connected graph, $\bar{t}_f$ and $\bar{\mathbf{x}}$ are broadcast to all agents.  Agents compute and implement the control inputs to reach $\bar{\mathbf{x}}$ at $\bar{t}_f$. 

 \begin{remark}
     The computation of $\bar{t}_f$ can be implemented in a distributed manner along similar lines as discussed in our previous article \cite{tacsubmitted2024}.
 \end{remark}
 
\begin{figure*}[h!]
    \begin{subfigure}{0.32\textwidth}
        \centering          \includegraphics[width=\linewidth]{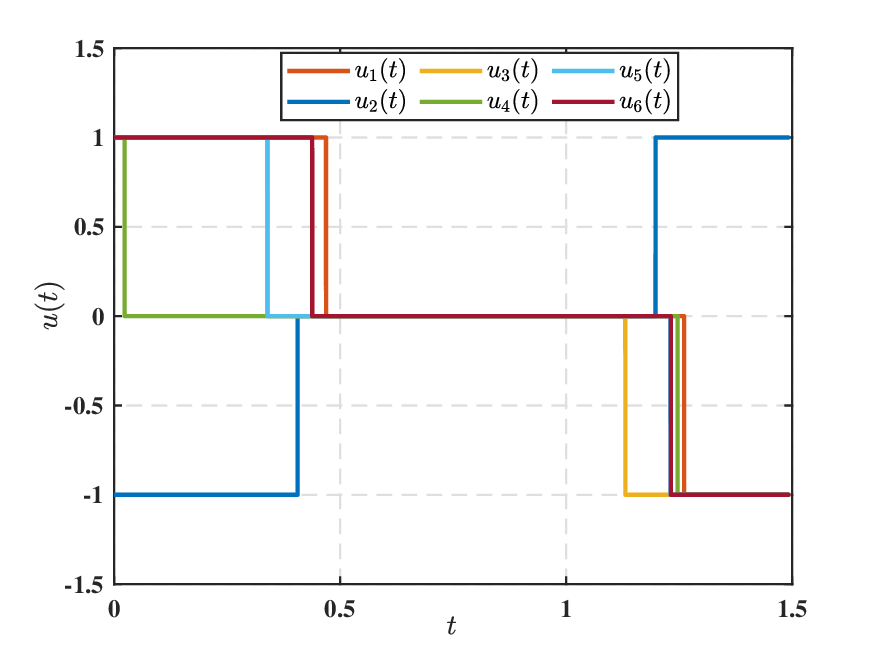}
        \caption{Control profile}
        \label{fig:3a}
    \end{subfigure}
    \hfill
    \begin{subfigure}{0.32\textwidth}
        \centering         \includegraphics[width=\linewidth]{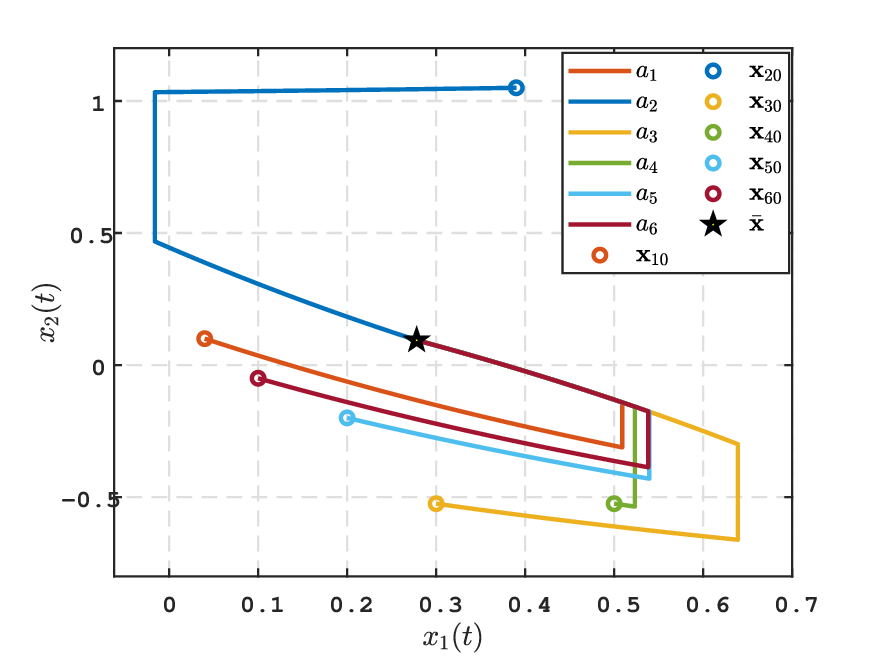}
        \caption{Phase plots}
        \label{fig:3b}
    \end{subfigure}
    \hfill
    \begin{subfigure}{0.32\textwidth}
        \centering             \includegraphics[width=\linewidth]{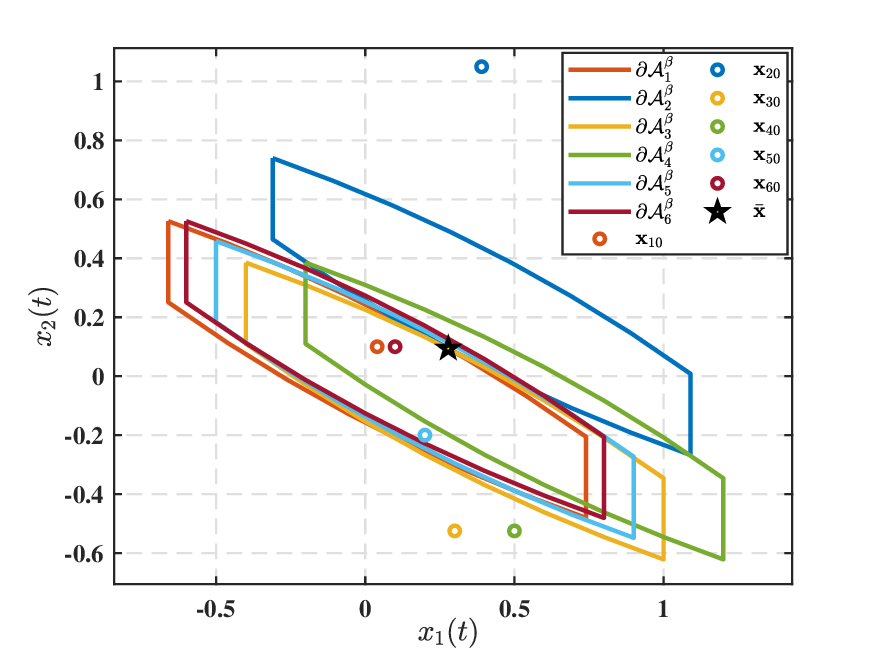}
       \caption{$\cap_{i=1}^{6}\partial {\mathcal{A}^{\beta}_i}$ with $\bar{\mathbf{x}} = [0.2781,0.0941]^{\top}$.}
        \label{fig:3c}
    \end{subfigure}
    \caption{State-trajectories, control profile and attainable set intersection for $N=6$ agents}
    \label{Fig:3}
\end{figure*}
\section{Example}
Consider a set of initial conditions for six agents, $a_1$, $a_2$, $a_3$, $a_4$, $a_5$, and $a_6$ i.e.,  $\mathbf{x}_{10} =[0.04,0.1]^{\top}$, $\mathbf{x}_{20} =[0.39,1.05]^{\top}$, $\mathbf{x}_{30} =[0.3,-0.525]^{\top}$, $\mathbf{x}_{40} =[0.5,-0.525]^{\top}$, $\mathbf{x}_{50} =[0.2, -0.2]^{\top}$, $\mathbf{x}_{60} =[0.1,-0.05]^{\top}$, respectively and $\beta = 0.7$. Note that ${x_1}_{\max}=0.5$ and ${x_1}_{\min}=0.04$. Thus, the set of initial conditions satisfies the inequality \eqref{eq:consensuspossible}. Among all $20$ triplets, the minimum consensus time is given by $\{a_1,a_2,a_3\}$, which is the maximum among all the other triplets of agents. Moreover, for $\{a_1,a_2,a_3\}$ triplet, the consensus point satisfies $\bar{\mathbf{x}} \in \partial \mathcal{A}_i^\beta(\bar{t}_f,\mathbf{x}_{0,i})$ whereas for remaining agents $\bar{\mathbf{x}} \in \text{int} \mathcal{A}_i^\beta(\bar{t}_f,\mathbf{x}_{0,i})$. Now, to ensure all agents reach the consensus point simultaneously at $\bar{t}_f$, the fuel budget $\beta$ for remaining agent is recomputed. This is done by substituting the computed values of $\bar{t}_f$ and $\bar{\mathbf{x}}$ in $\Gamma^{s_1}_{(L_i,x_{20,i})}(w_{1},x_{2},q_f)$, or  $\Gamma^{s_3}_{(L_i,x_{20,i})}(w_{1},x_{2},q_f)$ depending on the satisfaction of the associated inequality constraints. Subsequently, the corresponding switching instants are determined based on the recomputed $\beta$ using equations \eqref{eq:sw1},\eqref{eq:sw2}, \eqref{eq:sw3},\eqref{eq:sw4}, so that they reach the consensus point $\bar{\mathbf{x}}$ exactly at $\bar{t}_f$. The control profiles, the phase plots, and the attainable sets of respective agents are shown in Fig. \ref{fig:3a}, Fig. \ref{fig:3b} and Fig. \ref{fig:3c}, respectively. The intersection of the attainable sets of all agents indicates the minimum time consensus point which is $\bar{\mathbf{x}}=[ 0.2781, 0.0941]^{\top}\in \mathcal{X}_c^{\beta}(\infty)$. The minimum time to reach consensus $\bar{t}_f$ is $1.4918$. The control inputs and corresponding fuel utilization for all agents are listed in Table \ref{tab:fuel_utilization}. 
\begin{table}[ht]
\centering
\caption{\small Fuel utilization and control profiles of $N=6$ agents to reach $\bar{\mathbf{x}} = [0.2781,\, 0.0941]^{\top}$ at the minimum time $\bar{t}_f = 1.4918$.}
\label{tab:fuel_utilization}
\renewcommand{\arraystretch}{1.2} 
\setlength{\tabcolsep}{3.5pt}       
\begin{tabular}{@{} c c c c @{}}
\toprule
\textbf{Agent} & \textbf{Fuel Used ($\beta$)} & \textbf{Switching Times $(t_1, t_2)$} & \textbf{Control Profile $u(t)$} \\ 
\midrule
$a_1$ & 0.7000 & (0.4690, 1.2609) & $\{+1,\, 0,\, -1\}$ \\ 
$a_2$ & 0.7000 & (0.4060, 1.1978) & $\{-1,\, 0,\, +1\}$ \\ 
$a_3$ & 0.7000 & (0.3390, 1.1309) & $\{+1,\, 0,\, -1\}$ \\ 
$a_4$ & 0.2687 & (0.0234, 1.2465) & $\{+1,\, 0,\, -1\}$ \\ 
$a_5$ & 0.6009 & (0.3395, 1.2304) & $\{+1,\, 0,\, -1\}$ \\ 
$a_6$ & 0.6983 & (0.4382, 1.2317) & $\{+1,\, 0,\, -1\}$ \\ 
\bottomrule
\end{tabular}
\end{table}

\section{Conclusion}
For a set of $N$ second-order LTI agents, a method to compute the minimum time to consensus and the corresponding consensus point subject to a fixed fuel budget and bounded inputs is presented.  The key tool of our method, i.e., the attainable set for each agent under fuel budget constraints, is characterized. The convexity property of the attainable set allows us to use Helly's theorem for distributing the computation required in our method.  Initial conditions configurations for which consensus is possible under bounded inputs and fixed fuel budget constraints are also characterized. An extension to a general class of $n$-order LTI systems is currently under investigation and will be reported in the future.

\bibliographystyle{IEEEtran}
\bibliography{reference}

\end{document}